\documentclass[12pt]{article}
\usepackage{amsmath}

\usepackage[english]{babel}
\usepackage{amsmath}
\usepackage{amsthm,amssymb,latexsym}

\newtheorem{prop}{Proposition}
\newtheorem{lemma}[prop]{Lemma}

\newtheorem{theorem}[prop]{Theorem}
\theoremstyle{definition}

\title {There are no iterated morphisms that define the Arshon sequence and the $\sigma$-sequence}
\author {Sergey Kitaev}
\begin{document}
\maketitle

\begin{abstract} 
In \cite{3}, Berstel proved that the Arshon sequence cannot be obtained by iteration of a morphism. An alternative proof of this fact is given here. The $\sigma$-sequence was constructed by Evdokimov in order to construct chains of maximal length in the $n$-dimensional unit cube. It turns out that the $\sigma$-sequence has a close connection to the Dragon curve \cite{8}. We prove that the $\sigma$-sequence can not be defined by iteration of a morphism.
\end{abstract}

\section{Introduction and Background}

In 1937, Arshon gave a construction of a sequence of symbols $w$ over the alphabet $\{ 1,2,3 \}$, constructed as follows: Let $w_1=1$. For $k \ge 1$, $w_{k+1}$ is obtained by replacing the letters of $w_k$ in odd positions thus: 
$$1 \rightarrow 123, \ 2 \rightarrow 231, \ 3 \rightarrow 312$$ 

and in even positions thus: 
$$1 \rightarrow 321, \ 2 \rightarrow 132, \ 3 \rightarrow 213.$$

Then 
$$w_2=123, \ \ w_3=123132312,$$
and each $w_i$ is the initial subword of $w_{i+1}$, so the infinite symbolic sequence \linebreak $w=\lim\limits_{ n \to \infty }w_n$ is well defined. It is called the {\em Arshon sequence}.

This method of constructing $w$ is called the {\em Arshon Method} ({\bf AM}), and $\psi$ will denote the indicated map of the letters 1, 2, 3 according to position as described above. 

We will denote the natural decomposition of $w$ in $3$-blocks by lower braces:
$$w=\underbrace{123}\underbrace{132}\underbrace{312}\ldots.$$

The paper by Arshon \cite{2} was published in connection with the problem of constructing a nonrepetitive sequence on a 3-letter alphabet, that is, a sequence that does not contain any subwords of the type $XX=X^2$, where $X$ is any non-empty word of a 3-letter alphabet. The sequence $w$ has that property. The question of the existence of such a sequence, as well as the questions of the existence of sequences avoiding repetitions in other meanings, were studied in algebra \cite{1,9,10}, discrete analysis \cite{5,7} and in dynamical systems \cite{14}.

\

Any natural number $n$ can be presented  unambiguously as $n = 2^t(4s + \sigma)$, where $\sigma < 4$, and $t$ is the greatest natural number such that $2^t$ divides $n$. If $n$ runs through the natural numbers then $\sigma$ runs through the sequence that we will call  the {\em $\sigma$-sequence}. We let $w_{\sigma}$ denote that sequence. Obviously, $w_{\sigma}$ consists of $1$s and $3$s. The initial letters of $w_{\sigma}$ are $11311331113313 \ldots$. 

In \cite{6,17}, Evdokimov constructed chains of maximal length in the $n$-dimensional unit cube using the $\sigma$-sequence. Originally, the $\sigma$-sequence was defined by the following inductive scheme:

\centerline{$C_1=1$, \ \ \ \ $D_1=3$}

\centerline{$C_{k+1}=C_k1D_k$, \ \ \ \ $D_{k+1}=C_k3D_k$}

\centerline{$k=1, 2, \ldots$}

and $w_{\sigma}=\lim\limits_{ k \to \infty }{C_k}$.

Our definition above of the $\sigma$-sequence is equivalent to this one.

The motivation for studying the $\sigma$-sequence is that this sequence has a close connection to the well-known {\em Dragon curve}, discovered by physicist John E. Heighway and defined as follows: we fold a sheet of paper in half, then fold in half again, and again, etc. and then unfold in such way that each crease created by the folding process is opened out into a 90-degree angle. The ``curve'' refers to the shape of the partially unfolded paper as seen edge on. If one travels along the curve, some of the creases will represent turns to the left and others turns to the right. Now if 1 indicates a turn to the right, and 3 to the left, and we start travelling along the curve indicating the turns, we get the $\sigma$-sequence~\cite{8}.  

\  

Let $\Sigma$ be an alphabet and ${\Sigma}^{\star}$ be the set of all words of $\Sigma$. A map $\varphi: {\Sigma}^{\star} \rightarrow {\Sigma}^{\star}$ is called a {\em morphism}, if we have $\varphi(uv)=\varphi(u)\varphi(v)$ for any $u,v \in {\Sigma}^{\star}$. It easy to see that a morphism $\varphi$ can be defined by defining $\varphi(i)$ for each $i \in \Sigma$. The set of all rules $i \rightarrow \varphi(i)$ is called a {\em substitution system}. We create words by starting with a letter from the alphabet $\Sigma$ and iterating the substitution system. Such a substitution system is called a {\em D0L (Deterministic, with no context Lindenmayer) system} \cite{12}. D0L systems are a classical object of Formal Language Theory. They are interesting from mathematical point of view, but also have applications in theoretical biology~\cite{11}. 

Suppose a word $\varphi(a)$ begins with $a$ for some $a \in \Sigma$, and that the length of ${\varphi^k}(a)$ increases without bound. The symbolic sequence $\lim\limits_{ k \to \infty }{\varphi^k}(a)$ is called a {\em fixed point} of the morphism $\varphi$. 

We now study classes of sequences that are defined by iterative schemes. There are many techniques to study sequences generated by morphisms \cite{13}. It is reasonable to try to determine if a sequence under consideration can be obtained by iteration of a morphism.

Since the construction of the Arshon sequence $w$ is similar to the iterated morphism scheme, and because $w$ is constructed by two morphisms $f_1$ and $f_2$, applied depending on whether the letter position is even or odd, we might expect that there exists a morphism $f$ which generates $w$. But this turns out not to be true. In \cite{3}, Berstel proved that the Arshon sequence cannot be obtained by iteration of a morphism. However, Berstel's proof appeals to advanced machinery, namely a sort of Wilf-Fine theorem due to Cobham \cite{4}. In Theorem 1, we show that such machinery is not needed here, and offer a simple and direct proof. 

Naturally a question arises as to the possibility of constructing $w_{\sigma}$ using the iteration of a morphism, since such a construction could help us in studying $w_{\sigma}$, and thus in studying the Dragon curve. This also turns out not to be true, as shown in Theorem 10.
 
\section{The Arshon Sequence}

\begin{theorem} There does not exist a morphism whose fixed point is the Arshon sequence.\end{theorem}

{\bf Note.} A corollary of this theorem is the non-existence of a morphism whose iteration gives the Arshon sequence. In fact, if such a morphism exists, it must have the property that $1$ is mapped to $1X$ by the action of the morphism, where $X$ is some word, and from this it follows that the Arshon sequence is a fixed point of this morphism.

\begin{proof} It is enough to prove the non-existence of a morphism $f$ with the property $w=f(w)$, since from the definition of a fixed point we have that if $w$ is a fixed point of the morphism $f$ then $w = f(w)$. Suppose there exists a morphism $f$ such that $f(1)=X,\ f(2)=Y,\ f(3)=Z$ and $w = f(w)$. Moreover, in order to avoid dealing with the identity morphism, we assume that $|f(123)|>3$. We do this since $w$ is a fixed point of the identity morphism, but iterations of this morphism obviously do not give us the Arshon sequence. From all such morphisms $f$ we choose a morphism with minimal length of~$X$.

Clearly not all of $|X|$, $|Y|$ and $|Z|$ can be 0. Suppose for example then that $|X|=0$, $|YZ| \ne 0$. Then $w$ contains the subword $f(123132)=XYZXZY=YZZY$. If $|Z| \ne 0$ then $w$ contains the repetition $ZZ$. Otherwise $w$ contains repetition $YY$. This is impossible, hence $|X|\ne 0$. One can check in the same way that $|Y| \ne 0$ and $|Z| \ne 0$. Now $|X|+|Y|+|Z| \ne 3$, since otherwise $|f^{\ell}(1)|=1$ for ${\ell}=1,2...$, and $w$ is not a fixed point of the morphism $f$. Now
$$f(w) = w = XYZXZYZXY \ldots,$$
hence $X$ consists of $|X|$ of the first letters of $w$, $Y$ is $|Y|$ of the following letters, and $Z$ is $|Z|$ of the letters following that. 

We will use upper braces to show the decomposition of $w$ into $f$-blocks (that is, to show the disposition of the words $X$, $Y$ and $Z$ in $w$). We have 
$$ w = \overbrace{ \underbrace{123} \underbrace{132} \ldots a_{|X|}}^{X} \overbrace{ a_{|X|+1} \ldots a_{|X|+|Y|}}^{Y} \overbrace{ a_{|X|+|Y|+1} \ldots a_{|X|+|Y|+|Z|}}^{Z} \overbrace{ a_{|X|+|Y|+|Z|+1}\ldots }^{X} \ldots,$$ where all $a_i$ are letters of the alphabet $\{ 1, 2, 3 \}$.

\begin{lemma} We have $|X| + |Y| + |Z| \equiv 0 \pmod{3}$. \end{lemma}

\begin{proof} From the structure of $w$, the frequencies of 1, 2, 3 in $w$ coincide, hence the frequencies of these letters in $f(w) = w$ coincide as well. But this is only possible when $|X| + |Y| + |Z| \equiv 0 \pmod{3}$. Indeed, otherwise there are two letters, whose frequencies in $f(123)=XYZ$ do not coincide, which implies that the frequencies of these letters in $f(w) = w$ do not coincide as well, since $w$ can be written as $w=W_1W_2W_3 \ldots$, where $W_i$ is a permutation of the letters $X$, $Y$ and $Z$. \end{proof}

\begin{lemma} The situation $|X| \equiv |Y| \equiv |Z| \equiv 0 \pmod{3}$ is impossible.\end{lemma}

\begin{proof} Suppose $|X| \equiv |Y| \equiv |Z| \equiv 0 \pmod{3}$. Then $X$, $Y$ and $Z$ consist of a whole number of 3-blocks. The properties of $\psi$ mean that the morphism $g$ given by $g(1) = \psi^{-1}(X)$, $g(2) = \psi^{-1}(Y)$, $g(3) = \psi^{-1}(Z)$ satisfies $w = g(w)$. By the minimality of $f$, $|g(123)|=3$, so that $g(123)=123$. This implies that $f(1)=123$, $f(2)=132$, $f(3)=312$. However, this impossible, since then $f(1231)=123132312123$ is not a prefix $w$. \end{proof}

We define the {\em $N$th $f$-block $X$} to be the block $X$ obtained by considering the decomposition of $w$ into $f$-blocks and then taking the $N$th block that corresponds to the word $X$ (possibly skipping some blocks that corresponds to the words $Y$ and $Z$ without counting them). Thus, for example, when we say ``the 4th $f$-block $X$'', we mean not the following marked block $w = XYZ{\bf X}ZY \ldots$, but the following marked block $w = XYZXZYZXYZY{\bf X} \ldots$.

\begin{lemma} With the assumption of the existence of the morphism $f$, $|X|~\le~5$. \end{lemma}

\begin{proof} Suppose $|X| \ge 6$, that is, $X = 123132 \ldots $. If $|X| \equiv 2 \pmod{3}$ ($|X| \equiv 1 \pmod{3}$), then $|X| \ge 7$ and using Lemma 2 we consider the 4th $f$-block $X = 12\underbrace{313}\ldots$ ($X = 1\underbrace{231}\underbrace{323}\ldots$). This contradicts the {\bf AM}. Hence $|X| \equiv 0 \pmod{3}$. 

It follows from Lemma 3 that the situation $|Y| \equiv 0 \pmod{3}$ is impossible. If $|Y| \equiv 1 \pmod{3}$ ($|Y| \equiv 2 \pmod{3}$), then we consider the 10th (3rd) $f$-block $X = 12\underbrace{313}2\ldots$ and it brings us to a contradiction with the {\bf AM}. Hence if $|X| \ge 6$ then the morphism $f$ can not exist.  \end{proof}

\begin{lemma} With the assumption of the existence of the morphism $f$, $|X|~\ne~1$.\end{lemma}

\begin{proof} If $|X| = 1$, then $X = 1$ and the length of the words ${f^k}(1)$ for $k = 1, 2, \ldots$ does not increase, whence $w$ is not a fixed point of the morphism $f$. This is a contradiction.  \end{proof}

\begin{lemma} With the assumption of the existence of the morphism $f$, $|X|~\ne~2$.
\end{lemma}

\begin{proof} Suppose $|X| = 2$, that is $X = 12$.

We have $|X| \equiv 2 \pmod{3}$, hence, using Lemma 2, we have $|Y|+|Z| \equiv 1 \pmod{3}$. 

We consider the 2nd $f$-block $X$ and the $f$-block $Z$ next after it. It can be seen that $Z$ begins with 3. We consider the 4th $f$-block $X$ and $Y$ preceding it and find that $Y$ ends with 3. But then, considering $YZ$, which is a subword of $w$, we see, that 33 is a subword of $w$, which is impossible. That is for $|X| = 2$ the morphism $f$ cannot exist. \end{proof}

The 3-blocks 123, 231, 312 are said to be {\em odd} 3-blocks. All other 3-blocks are said to be {\em even}.

\begin{lemma} With the assumption of the existence of the morphism $f$, $|X|~\ne~3$.\end{lemma}

\begin{proof} Suppose $|X| = 3$, that is $X = 123$.

We have $|X| \equiv 0 \pmod{3}$, hence, using Lemma 2 we have $|Y|+|Z| \equiv 0 \pmod{3}$. Considering the {\bf AM}, the 2nd $f$-block $X$ must be an odd 3-block, hence $|Y|+|Z| \equiv 1 \pmod{2}$. 

Let $|Z| \ge 2$. Then the 2nd $f$-block $Z$ begins with an even 3-block, and the 3rd $Z$ begins with an odd 3-block. This is impossible. (Note that 2 letters define the evenness of the 3-block unambiguously.) Thus $|Z| = 1$.

Let $|Y| \ge 2$. In $XYZ$ (or in an arbitrary permutation of these letters) there is an even number of 3-blocks, so the 9th $f$-block $Y$ begins with an odd 3-block, but the 1st $Y$ begins with an even 3-block. Hence $|Y| = 1$.

This is a contradiction with $|Y|+|Z| \equiv 0 \pmod{3}$ (and also a contradiction with $|Y|+|Z| \equiv 1 \pmod{2}$). That is for $|X| = 3$ the morphism $f$ cannot exist.  \end{proof}

\begin{lemma} With the assumption of the existence of the morphism $f$, $|X|~\ne~4$.\end{lemma}

\begin{proof} Suppose $|X| = 4$, that is $X = 1231$.

We have $|X| \equiv 1 \pmod{3}$, hence, using Lemma 2, we have $|Y|+|Z| \equiv 2 \pmod{3}$. 

We have $|Y| \ge 2$, since otherwise $Y = 3$ and hence $XYX$ which is a subword of $w$, contains 3131, which is impossible. Hence $Y = 32 \ldots $. We consider $ZX$ and $ZY$ and see that $Z$ ends with 2, since otherwise $w$ contains $Z1$ and $Z3$, and thus the square 11 or 33. Now $|Z| \ge 2$, since otherwise $Z = 2$ and $XZX$ which is a subword of $w$, contains 1212, which is impossible. Hence $Z = \ldots 32$, or $Z = \ldots 12$. The former is impossible since 3232 is contained in $ZY$, and hence in $w$. The latter is impossible too, since considering the 9th $f$-block $Z$ and the $f$-block $X$ following it, we obtain $ZX = \ldots \underbrace{121}\underbrace{231}$, which contradicts the {\bf AM}. That is for $|X| = 4$ the morphism $f$ cannot exist.  \end{proof}

\begin{lemma} With the assumption of the existence of the morphism $f$, $|X|~\ne~5$. \end{lemma}

\begin{proof} Suppose $|X| = 5$, that is $X = 12313$.

We have $|X| \equiv 2 \pmod{3}$, hence, using Lemma 2, we have $|Y|+|Z| \equiv 1 \pmod{3}$. Then the 4th $f$-block is $X = 12\underbrace{313}$, which is a contradiction with the {\bf AM}. That is if $|X| = 5$ then the morphism $f$ cannot exist.  \end{proof}

From Lemmas 4 -- 9 we have a contradiction with the assumption of the existence of the morphism $f$. This proves Theorem 1. \end{proof}

{\bf Remark.} In \cite{2}, Arshon gave the construction of a nonrepetitive sequence $w_n$ for an $n$-letter alphabet, where $n$ is any natural number greater than or equal to 3. It is easy to see that, for even $n$, there exists a morphism $f_n$ that defines $w_n$. Namely, for $1 \leq i \leq n$, one has:

\[ f_n(i) = \left\{ \begin{array}{ll}
i(i+1) \ldots n12 \ldots i-1, & \mbox{if $i$ is odd,} \\
(i-1)(i-2) \ldots 1n(n-1) \ldots i, & \mbox{if $i$ is even.}
\end{array}
\right. \] 

Theorem 1 shows that for $n=3$ such a morphism does not exist. However, whether there exists a morphism defining $w_n$ for arbitrary odd $n$ is still an open question.
   
\section{The $\sigma$-sequence}

\begin{theorem} There does not exist a morphism whose iteration defines the sequence $w_{\sigma}$. \end{theorem}

\begin{proof} Suppose there exists a morphism $f$, such that $f(1)=X$, $f(3)=Y$ and $w_{\sigma} = \lim\limits_{ k \to \infty }{f^k}(1)$. Obviously, $X$ consists of the first $|X|$ letters of $w$, where $|X|$ is the length of $X$. 

\begin{lemma} The subsequence of $w_{\sigma}$ consisting of the letters in odd positions is the alternating sequence of $1$s and $3$s: $1313131 \ldots$. \end{lemma}

\begin{proof} The odd positions of $w_{\sigma}$ correspond to the odd numbers $n=2^0(4s + \sigma)=4s+\sigma$, so clearly $\sigma$ alternates between $1$ and $3$. \end{proof}

\begin{lemma} If there exists a morphism $f$ whose iteration gives $w_{\sigma}$ then $|X| \equiv 0 \pmod{4}$. \end{lemma}

\begin{proof} It is easy to see that $f(1) = 1X^{(1)}$, where $|X^{(1)}| \ge 1$, since otherwise $|f^k(1)| = 1$, for $k = 1, 2, 3 \ldots $, so $w_{\sigma}$ cannot be obtained by iterating $f$.

Suppose $|X^{(1)}| = 1$, that is $f(1) = 11$. But then $w_{\sigma}$ consists of $1$s only, which is impossible, hence $f(1) = 11X^{(2)}$, where $|X^{(2)}| \ge 1$.

Suppose $|X^{(2)}| = 1$, that is $f(1) = 113$. Since $w_{\sigma}$ has the subword $111$, then $w_{\sigma}$ has a subword $f(111) = 113113113$. If $f(111)$ begins with a letter in an odd position, then the marked letters ${\bf 1}1{\bf 3}1{\bf 1}3{\bf 1}1{\bf 3}$, read from left to right will make up consecutive letters of $w_{\sigma}$ in odd positions. This contradicts Lemma 11. If $f(111)$ begins with a letter in an even position, then marking letters in odd positions will lead to the same contradiction with Lemma 11, hence $f(1) = 113X^{(3)}$, where $|X^{(3)}| \ge 1$.

Suppose $|X^{(3)}| = 1$, that is $f(1) = 1131$. Then $f^2(1) = 11311{\bf 1}31Y1131$ and the marked letter does not coincide with the letter of $w_{\sigma}$ standing in the same place, hence $f(1) = 1131X^{(4)}$, where $|X^{(4)}| \ge 1$. 

If $|X|$ is odd, then the marked letters in $f^2(1) = 1131X^{(4)}1{\bf 1}3{\bf 1}X^{(4)} \ldots$ are two consecutive letters in odd places. This contradicts Lemma 11. Hence $|X|$ is even.

We have $f^2(1) = XX \ldots = X1131X^{(4)} \ldots$, whence the next-to-last letter of $X$ is in an odd position and is equal to $3$, since otherwise two consecutive $1$ in $w_{\sigma}$ stand at odd places, which contradicts Lemma 11. The natural number which corresponds to the next-to-last letter of $X$ is written as $2^0(4s + 3)$, the next number is equal to $|X|$ and to $2^0(4s + 3) + 1 = 4(s + 1) \equiv 0 \pmod{4}$.  \end{proof}

The following Lemma is straightforward to prove.

\begin{lemma} If $n_1 = 2^{t_1}(4s_1 + 1)$, $n_2 = 2^{t_2}(4s_2 + 1)$, $n_3 = 2^{t_3}(4s_3 + 3)$ and $n_4 = 2^{t_4}(4s_4 + 3)$ then $n_1n_2$, $n_3n_4$ can be written as $2^t(4s + 1)$, and $n_1n_3$ as $2^t(4s + 3)$. \end{lemma}

It follows from Lemma 12 that $|X| = 4t$. 

Suppose $X$ ends with $1$ (the case when $X$ ends with $3$ is similar), that is at the $(4t)$th position in $X$ we have $1$. According to the multiplication by 2 does not change $\sigma$, so at the $(2t)$th position in $X$ we have $1$. 

Consider $f^2(1) = X{\bf X} \ldots $. The letters of the marked $X$ occupy the positions of $f^2(1)$ from $(4t+1)$th to $(8t)$th. Since $X={\bf X}$, then at the $(6t)$th place we have $1$. But $6t = 3(2t)$, whence, by Lemma 13, at the $(2t)$th and the $(6t)$th places there must stand different letters. This is a contradiction and Theorem 10 is proved. \end{proof}

\end{document}